\documentclass[graybox]{svmult}
\usepackage{amsmath, amssymb, amsopn}
\usepackage{array}
\usepackage{amscd}
\usepackage{epsfig}
\usepackage{color}
\usepackage{comment}
\usepackage{setspace}
\usepackage{multirow}
\usepackage{booktabs}
\usepackage{todonotes}
\usepackage{mathtools}
\usepackage{url}

\usepackage{type1cm}        
\usepackage{makeidx}         
\usepackage{graphicx}        
\usepackage{multicol}        
\usepackage[bottom]{footmisc}

\usepackage{newtxtext}       %
\usepackage{newtxmath}       

\makeindex             

\DeclareMathOperator{\Div}{div}

\DeclareMathOperator{\Grad}{\nabla}

\DeclareMathOperator{\ssum}{\textstyle \sum}
\DeclareMathOperator{\diag}{diag}

\newcommand{\inner}[2]{\left( #1, #2 \right)}
\newcommand{\foralls}{\forall \,}

\newcommand{\R}{{\mathbb R}}

\newcommand{\exchange}{{E}}
\newcommand{\foo}{{K}}
\allowdisplaybreaks
\begin{document}

\title*{Parameter robust preconditioning for multi-compartmental Darcy equations}
\author{Eleonora Piersanti, Marie E. Rognes, and Kent-Andre Mardal}
\authorrunning{E. Piersanti et al.}

\institute{E. Piersanti \at Simula Research Laboratory, 1325 Lysaker, Norway  \email{eleonora@simula.no}
\and M. E. Rognes \at Simula Research Laboratory, 1325 Lysaker, Norway \email{meg@simula.no}
\and K.-A. Mardal \at Department of Mathematics, University of Oslo, 0316 Oslo, Norway \email{kent-and@math.uio.no}}
%
%
\maketitle

\abstract{In this paper, we propose a new finite element solution
  approach to the multi-compartmental Darcy equations describing flow
  and interactions in a porous medium with multiple fluid
  compartments. We introduce a new numerical formulation and a
  block-diagonal preconditioner. The robustness with respect to
  variations in material parameters is demonstrated by
  theoretical considerations and numerical examples.}

\section{Introduction}
The multi-compartment Darcy equations\footnote{In this paper, we will
  also refer to these equations as the multiple--network porosity
  (MPT) equations.} extend the single compartment Darcy model and
describe fluid pressures in a rigid porous medium permeated by
multiple interacting fluid networks. These equations have been used to
model perfusion in e.g.~the heart~\cite{MichlerEtAl2013,
  LeeEtAlDarcy2015}, the brain~\cite{JoszaEtAl2019} and the
liver~\cite{BrasnovaEtAl2018}. 
The static variant of the equations read as follows: for a
given number of networks $J \in \mathbb{N}$, find the network
pressures $p_{j}$ for $j = 1, \dots, J$ such that
\begin{equation}
  \label{eq:mpt}
    - K_j \Div \Grad p_{j} + \sum_{i=1}^J \xi_{j \leftarrow i} (p_j - p_i) = g_{j} \quad \text{in } \Omega,
\end{equation}
where $p_{j} = p_{j}(x)$ for $x \in \Omega \subset \R^{d}$ ($d = 1, 2,
3$), and $\Omega$ is the physical domain. The scalar parameter $K_j >
0$ represents the permeability of each network $j$. The parameter
$\xi_{j \leftarrow i} \geq 0$ is the exchange coefficient into network
$j$ from network $i$. These are assumed to be symmetric: $\xi_{j
  \leftarrow i} = \xi_{i \leftarrow j}$. The right hand side $g_j$ can
be interpreted as a source/sink term for each $j$. For simplicity, let
$p_j = 0$ on $\partial \Omega$ for $1 \leq j \leq J$.

The system of equations is elliptic as long as $K_j > 0$, but for $K_j
\ll \xi_j$ the diagonal dominance is lost for smooth components for
which $\|K_j^{1/2} \Grad p_{j}\| \le \|\xi_j^{1/2} p_j\|$. As diagonal
dominance is often exploited in standard preconditioning algorithms
such as for example multigrid, the consequence is a loss of
performance.  Here, we will therefore propose a transformation of the
system of equations that enable the use of standard preconditioners.
In detail, we propose and analyze a new approach to constructing
finite element formulations and associated block--diagonal
preconditioners of the system~\eqref{eq:mpt}. The key idea is to
change variables through a transformation $T$ that gives simultaneous
diagonalization by congruence of the operators involved. We preface
and motivate the new approach by a demonstration of lack of robustness
of a standard formulation for high exchange parameters.

\section{Lack of parameter robustness in standard formulation}
\label{sec:motivation}

A standard variational formulation of~\eqref{eq:mpt} reads as follows:
find $p_j \in H^1_0 = H^1_0(\Omega)$ for $1 \leq j \leq J$ such that:
\begin{equation}
  \label{eq:mpt:vf}
  \inner{ K_{j} \Grad p_j}{\Grad q_j}
  + \ssum_{i=1}^J \inner{\xi_{j \leftarrow i} (p_j - p_i)}{q_j}
  = \inner{g_j}{q_j} \quad \foralls q_j \in H^1_0,
\end{equation}
where $\inner{\cdot}{\cdot}$ denotes the $L^2(\Omega)$ inner product.
The system~\eqref{eq:mpt:vf} can be written in the alternative form:
\begin{equation}
  \label{eq:mpt:vf2}
  k (\bold{p}, \bold{q}) + e(\bold{p}, \bold{q}) = \inner{\bold{g}}{\bold{q}},
\end{equation}
with $\bold{p} = (p_1, p_2, \dots, p_J),$ $\bold{q} = (q_1, q_2,
\dots, q_J),$ $\bold{g} = (g_1, g_2, \dots, g_J)$, and with matrix
form
\begin{equation*}
  \mathcal{A} \bold{p} = \bold{g},
\end{equation*}
where 
\begin{equation*}
\label{eq:mpt:KD}
\mathcal{A} = \mathcal{K} +  \exchange = \begin{pmatrix} K_1 \Delta & 0 & \cdots & 0 \\
0 & K_2 \Delta  & \cdots & 0 \\
\vdots & \vdots & \ddots & \vdots \\ 
0 & 0 & \cdots & K_J\Delta \\
\end{pmatrix} + \begin{pmatrix} \sum_{i=1}^J \xi_{1 \leftarrow i} & - \xi_{1 \leftarrow 2} & \cdots & - \xi_{1 \leftarrow J}  \\
- \xi_{1 \leftarrow 2} & \sum_{i=1}^J \xi_{2 \leftarrow i} & \cdots & - \xi_{2 \leftarrow J} \\
\vdots & \vdots & \ddots & \vdots \\ 
- \xi_{1 \leftarrow J} & - \xi_{2 \leftarrow J} & \cdots & \sum_{i=1}^J \xi_{J \leftarrow i} \\
\end{pmatrix}.
\end{equation*}

Taking the blocks on the diagonal of $\mathcal{A}$ we can immediately
define a block diagonal preconditioner $\mathcal{B}$:
\begin{equation}
  \label{eq:mpt:precond1}
  \mathcal{B} =
  \diag \left( -K_1\Delta + \sum_{i=1}^J \xi_{1 \leftarrow i},
  - K_2 \Delta + \sum_{i=1}^J \xi_{2 \leftarrow i}, \cdots,
  -K_J \Delta + \sum_{i=1}^J \xi_{J \leftarrow i} \right)
\end{equation}
Alas, this formulation and preconditioner is not robust for high
exchange parameters as illustrated by the following example.
\begin{example}
  \label{ex:mpt}
  In this example we illustrate the poor performance of the block
  diagonal preconditioner \eqref{eq:mpt:precond1} for the standard
  finite element discretization of the MPT equations~\eqref{eq:mpt}
  with $J=2$. In particular, we show that the proposed preconditioner
  is not robust with respect to the exchange coefficient $\xi_{1
    \leftarrow 2}$ and mesh refinement.  Let $\Omega = [0, 1]^2
  \subset \R^2$, and let $K_1 = K_2 = 1.0$, $g_j = 0$. To discretize
  the pressures $p_1, p_2$ we consider continuous piecewise linear
  finite elements defined relative to a $2 N \times N$ triangular mesh
  of $\Omega$. The results in Table \ref{tab:mpt} show that both the
  number of iterations and condition numbers increase somewhat less than linearly (predicted 
  by our theoretical analysis) in
  $\xi_{1 \leftarrow 2},$ for $\xi_{1 \leftarrow 2}$ above a threshold
  $>100$. The number of iterations also grow for increasing $N$
  (decreasing mesh size $h$) in this case.
\begin{table}[h]
\begin{center}
\begin{tabular}{l|lllll}
$\xi_{1 \leftarrow 2}$ &  \multicolumn{5}{c}{N}\\
\toprule
 & $8$ & $16$ & $32$ & $64$ & $128$\\
\midrule
 $10^{-2}$ & $3$ ($1.0$)      & $4$ ($1.1$)      & $4$ ($1.1$)       & $4$ ($1.1$)      & $4$ ($1.1$)       \\
 $10^0$ & $4$ ($1.1$)      & $4$ ($1.1$)       & $5$ ($1.1$)       & $5$ ($1.1$)      & $5$ ($1.1$)       \\
 $10^2$ & $29$ ($11$)    & $30$ ($11$)    & $28$ ($11$)     & $25$ ($11$)    & $24$ ($11$)     \\
 $10^4$ & $215$ ($1053$) & $740$ ($1026$) & $1131$ ($1012$) & $1232$ ($1014$) & $1058$ ($1014$) \\
 $10^6$  & $7$ ($2.0$)       & $20$ ($581$)   & $84$ ($686$)    & $394$ ($1140$) & $1467$ ($1755$) \\
\bottomrule
\end{tabular}
\caption{Number of iterations (and condition number estimates) of a CG
  solver of the system~\eqref{eq:mpt} with an algebraic multigrid
  (Hypre AMG) preconditioner of the form~\eqref{eq:mpt:precond1} with
  a random initial guess. Results for $\xi_{1 \leftarrow 2} = 10^{-4},
  10^{-6}$ are nearly identical to the $10^{-2}$ case.}
\label{tab:mpt}
\end{center}
\end{table}
\end{example}

We can examine Example~\eqref{ex:mpt} analytically. Define the induced norm
\begin{equation}
\label{eq:norm:Bmpt1}
\lVert \bold{p} \rVert^2_{\mathcal{B}} = \inner{\mathcal{B} \bold{p}}{\bold{p}}
= \ssum_{j=1}^J \inner{K_j \Grad p_j}{\Grad p_j} + \xi_j \inner{p_j}{p_j},
\end{equation}
where $\xi_j = \sum_{i=1}^J \xi_{j \leftarrow i}.$ We can show that
there exists an $\alpha > 0$ such that
\begin{equation}
\inner{\mathcal{A} \bold{p}}{\bold{p}} \geq \alpha \inner{\mathcal{B} \bold{p}}{\bold{p}}
\end{equation} 
for all $\bold{p}$, but depending on $K_j$ and $\xi_{j \leftarrow i}$,
as follows. Note that for all $\bold{p}$
\begin{equation}
\label{eq:MPT:coercivity1} 
\inner{\mathcal{A} \bold{p}}{\bold{p}}
= \inner{(\mathcal{K} + \exchange) \bold{p}}{\bold{p}}
\geq  \inner{\mathcal{K} \bold{p}}{\bold{p}},
\end{equation}
since 
\begin{equation*}
  \inner{\exchange \bold{p}}{\bold{p}}
  = \ssum_{i=1}^J \ssum_{j=1}^J  \inner{\xi_{j \leftarrow i}(p_j-p_i)}{p_j}
  = \frac{1}{2} \ssum_{j=1}^J \ssum_{i=1}^J \xi_{j\leftarrow i} \lVert p_j -p_i\rVert^2 \geq 0.
\end{equation*}
By definition and by applying the Poincar\'e inequality, we find that
there exists a constant $C_\Omega$ depending on the domain $\Omega$,
such that
\begin{equation}
\label{eq:MPT:coercivity2} 
\inner{\mathcal{K} \bold{p}}{\bold{p}}
= \sum_{j=1}^J \frac{K_j}{2} \lVert \Grad p_j \rVert^2 + \frac{K_j}{2} \lVert \Grad p_j \rVert^2 
\geq  \frac{1}{2} \sum_{j=1}^J K_j \lVert \Grad p_j \rVert^2 + \frac{C_\Omega K_j}{\xi_j} \xi_j \| p_j \|^2.
\end{equation}
Thus, using the definition of $\mathcal{B}$, we obtain that
\begin{equation}
  \inner{\mathcal{K} \bold{p}}{\bold{p}}
  \geq 
  \frac{1}{2} \min{\left(1, \min_j \frac{C_\Omega K_j}{\xi_j}\right) }
  \inner{\mathcal{B}\bold{p}}{\bold{p}}.
\end{equation}
We observe that the coercivity constant depends on the permeability
and exchange parameters and is such that it vanishes for vanishing
ratios of $K_j$ to $\xi_j.$ 

We can also show that there exists a constant $\beta$ such that
\begin{equation}
  \inner{\mathcal{A} \bold{p}}{\bold{q}} \leq \beta
  \|\bold{p}\|_{\mathcal{B}}\|\bold{q}\|_{\mathcal{B}} .
  \label{eq:continuity:1}
\end{equation} 
For any $\bold{p}$ and $\bold{q}$, applying the Cauchy--Schwartz
inequality twice we obtain
\begin{equation*}
  (\mathcal{A} \, \bold{p}, \bold{q})
  \leq \ssum_{j=1}^J\left( K_j \|\Grad p_j\| \|\Grad q_j\| + \ssum_{i=1}^J \xi_{j \leftarrow i}(\|p_j\|+\|p_i\|) \|q_j\|\right).
\end{equation*}
Applying the Cauchy--Schwartz inequality, the diffusion term is bounded as follows 
\begin{equation*}
  \ssum_{j=1}^J K_j \|\Grad p_j\| \|\Grad q_j\| \leq
  \left(\ssum_{j=1}^J K_j \|\Grad p_j\|^2\right)^{1/2} \left(\ssum_{j=1}^J K_j \|\Grad q_j\|^2\right)^{1/2} .
\end{equation*}
For the exchange terms, we can use the Cauchy-Schwartz inequality, the
symmetry of the exchange coefficients and Chebyshev's inequality to
show that
\begin{equation*}
  \ssum_{j=1}^J \ssum_{i=1}^J \xi_{j \leftarrow i} \|p_i\| \|q_j\| \leq
  J \left(\ssum_{j=1}^J \xi_j\|p_j\|^2\right)^{1/2} \left(\ssum_{j=1}^J \xi_j\|q_j\|^2\right)^{1/2},
\end{equation*}
and similarly for $\|p_j\|$ in place of $\| p_i
\|$. Thus~\eqref{eq:continuity:1} holds with continuity constant
$\beta$ equal to $J+1$.

The condition number of the preconditioned continuous system can be
estimated as the ratio between \eqref{eq:continuity:1} and
\eqref{eq:MPT:coercivity2}, c.f.~for example \cite{MardalWinther2011},
and tends to $\infty$ as $\xi_{j \leftarrow i} \rightarrow \infty$.
CG convergence is governed by the square root of the condition number
which in Example \ref{ex:mpt}, explains how the number of iterations
increase as $\xi_{1 \leftarrow 2}$ grows in Table~\ref{tab:mpt}.

\section{Change of variables yields parameter robust formulation}

In this section, we present a new approach to variational formulations
for the MPT equations. The key idea is to change from variables
$\bold{p}$ to variables $\tilde{\bold{p}}$ via a transformation $T$
such that the equation operators decouple. We can show that this is
always possible by simultaneous diagonalization of matrices by
congruence.

To this end, we define $\tilde{\bold{p}}$ and $\tilde{\bold{q}}$ as a
new set of variables such that
\begin{equation}
  \label{eq:pqtilde}
  \bold{p} = T \tilde{\bold{p}}, \quad \bold{q} = T \tilde{\bold{q}}.
\end{equation}
for a linear transformation map (matrix) $T : \R^J \rightarrow \R^J$
to be further specified. Substituting \eqref{eq:pqtilde} into
\eqref{eq:mpt:vf2}, we obtain a new variational formulation reading
as: find $\bold{\tilde{p}} \in (H^1_0)^J$ such that
\begin{equation}
\label{eq:mpt:vf:transformed}
k(T\bold{\tilde{p}}, T\bold{\tilde{q}}) + e(T\bold{\tilde{p}}, T\bold{\tilde{q}})
= \inner{ T^T\bold{g}}{\bold{\tilde{q}}} \quad \foralls \bold{\tilde{q}} \in (H^1_0)^J .
\end{equation}
The matrix form of the system is
\begin{equation}
\label{eq:mpt:coeffmatrix1:transformed}
\mathcal{\tilde{A}} \bold{\tilde{p}} = (\mathcal{\tilde{K}} + \tilde{E}) \bold{\tilde{p}}= T^T \bold{g} = \bold{\tilde{g}},
\end{equation}
where
\begin{equation}
  \label{eq:mpt:coeffmatrix2:transformed}
  \mathcal{\tilde{K}} = (-\Delta) \tilde{\foo}, \quad \tilde{\foo} = T^T \foo T, \quad 
  \tilde{\exchange}  = T^T \exchange T,
\end{equation}
where the matrix $\exchange \in \R^J \times \R^J$ is given in
Section~\ref{sec:motivation} and where we write $\foo = \diag (K_1,
K_2, \dots, K_J)$.

The key question is now whether there exists an (invertible)
transformation $T$ that simultaneously diagonalizes (by congruence)
$\foo$ and $\exchange$? More precisely, is there a matrix $T \in \R^J
\times \R^J$ such that
\begin{align}
  \label{eq:diagonalized}
  \tilde{K} = \diag(\tilde{K}_1, \tilde{K}_2, \dots, \tilde{K}_J), \quad
  \tilde{E} = \diag(\tilde{\xi}_1, \tilde{\xi}_2, \dots, \tilde{\xi}_J) \quad ?
\end{align}
By matrix analysis theory, see e.g.~\cite[Theorem 4.5.17,
  p.~287]{horn1990matrix}, there exists indeed such a $T$ since $K$ is
diagonal and non-singular and $E$ is symmetric and thus $C = K^{-1} E$
is diagonalizable. In particular, consider the case where $C$ has $J$
distinct eigenvalues $\lambda_j$ and eigenvectors $v_j$ for $j = 1,
\dots, J$. By taking $T = [v_1, v_2, \dots,
  v_J]$,~\eqref{eq:diagonalized} holds. Moreover, the eigenvalues
$\lambda_j$ are all real.

\begin{example}
  \label{ex:mpt:2net:diagonalization}
  To exemplify, we here show the diagonalization by congruence of a
  general 2--network system explicitly. Let
\begin{equation*}
  K = \begin{pmatrix} K_1 & 0 \\
    0 & K_2 
  \end{pmatrix}, \quad
  E = \begin{pmatrix} \xi_{1 \leftarrow 2} & - \xi_{1 \leftarrow 2} \\
    - \xi_{1 \leftarrow 2} & \xi_{1 \leftarrow 2}
  \end{pmatrix}.
\end{equation*}
Then,
\begin{equation*}
  C = K^{-1}E =
  \begin{pmatrix}
    \xi_{1 \leftarrow 2}/K_1 & -\xi_{1 \leftarrow 2}/K_1 \\
    -\xi_{1 \leftarrow 2}/K_2 &  \xi_{1 \leftarrow 2}/K_2, 
  \end{pmatrix}
\end{equation*}
has eigenvalues $e_1 = 0$ and $e_2 = \xi_{1 \leftarrow 2}(K_1 +
K_2)/(K_1 K_2)$ and the eigenvectors form the columns of $T$:
\begin{equation*}
  T =
  \begin{pmatrix}
    1 & K_2(\xi_{1 \leftarrow 2}/K_2 - \xi_{1 \leftarrow 2}(K_1 + K_2)/(K_1 K_2))/\xi_{1 \leftarrow 2} \\
    1 & 1
  \end{pmatrix},
\end{equation*}
Finally, we can verify that
\begin{align*}
  \tilde{K} & = T^T K T   
  = \begin{pmatrix} K_1 + K_2 & 0 \\
    0 & K_2 (K_1 + K_2)/K_1
  \end{pmatrix},\\ 
  \tilde{E} &= T^T E T =  \begin{pmatrix} 0 & 0 \\
    0 & \xi_{1 \leftarrow 2}(K_1^2 + K_1 K_2 + K_2(K_1 + K_2))/K_1^2
  \end{pmatrix} .
\end{align*}  
\end{example}

As the transformed system is diagonal and decoupled, a block--diagonal
preconditioner is readily available. In particular, we define
\begin{equation}
  \label{eq:mpt:transformed:precond}
  \mathcal{\tilde{B}} =  \mathcal{\tilde{A}} = \diag \left (
  -\tilde{K}_1\Delta + \tilde{\xi}_1,
  -\tilde{K}_2 \Delta + \tilde{\xi}_2,
  \dots,
  -\tilde{K}_J \Delta + \tilde{\xi}_J \right ).
\end{equation}
with norm
\begin{equation}
\lVert \bold{\tilde{p}} \rVert^2_{\mathcal{\tilde{B}}} = \inner{\tilde{\mathcal{B}} \bold{\tilde{p}}}{\bold{\tilde{p}}} = \sum_{j=1}^J \inner{\tilde{K}_j \Grad \tilde{p}_j}{\Grad \tilde{p}_j} + \tilde{\xi}_j \inner{\tilde{p}_j}{\tilde{p}_j}.
\end{equation}
Clearly, by definition, $\mathcal{\tilde{A}}$ and
$\mathcal{\tilde{B}}$ are trivially spectrally equivalent (with upper
and lower bounds independent of the material parameters).


\section{Numerical examples for the new formulation}

In this section, we present numerical results supporting the
theoretical considerations. All numerical experiments have been
conducted using a finite element discretization, using the FEniCS
library~\cite{LoggEtAl2012} and the \textrm{cbc.block}
package~\cite{MardalEtAl2012}. To discretize the pressures $p_j$ and
the transformed variables $\tilde{p}_j$, we consider continuous
piecewise linear ($P_1$) finite elements defined relative to each mesh
$\mathcal{T}_h$ of the domain $\Omega = [0, 1]^2$. We impose
homogeneous Dirichlet conditions on the whole boundary, and zero right
hand side(s). The linear systems were solved using a conjugate
gradient (CG) solver, with algebraic multigrid (Hypre AMG) with the
respective preconditioners, starting from a random initial guess.
The tolerance is set to $10^{-9},$ iterations are stopped at 3000, the condition number is just 
an estimation provided by the Krylov spaces involved in the iterations and will be lower than the real value.

\begin{example}
\label{ex:mpt:2net}
\begin{table}[t]
\centering
\begin{tabular}{l|l|l|l|l|l|l}
 $\xi_{1\leftarrow2}$ & $K_2$        & $N=8$            & $N=16$           & $N=32$           & $N=64$           & $N=128$          \\
\toprule
 \multirow{7}{*}{$10^4$}  & $10^{-6}$  & $277$ ($2139$)   & $1178$ ($2135$)   & $2395$ ($2035$)  & $3001$ ($2034$)  & $3001$ ($2034$) \\
   						  & $10^{-4}$  & $280$ ($2139$)   & $1180$ ($2135$)   & $2283$ ($2035$)  & $2860$ ($2034$)  & $3001$ ($2034$) \\
    					  & $10^{-2}$  & $275$ ($2117$)   & $1181$ ($2113$)   & $2325$ ($2014$)  & $2859$ ($2013$)  & $2988$ ($2011$) \\
    					  & $10^0$     & $242$ ($1054$)   & $935$ ($1026$)    & $1629$ ($1012$)  & $1556$ ($1014$)  & $1557$ ($1014$) \\
                          & $10^2$     & $62$ ($21$)      & $74$ ($22$)       & $74$ ($22$)      & $66$ ($22$)      & $64$ ($22$)     \\
                          & $10^4$     & $12$ ($1.6$)       & $11$ ($1.6$)        & $11$ ($1.6$)       & $11$ ($1.6$)       & $10$ ($1.6$)      \\
                          & $10^6$     & $7$ ($1.1$)        & $7$ ($1.1$)         & $7$ ($1.1$)        & $7$ ($1.1$)        & $7$ ($1.1$)       \\
\midrule
\multirow{7}{*}{$10^6$}  & $10^{-6}$ & $138$ ($34499$)  & $692$ ($42936$)   & $2999$ ($45584$) & $3001$ ($17128$) & $3001$ ($5730$) \\   
   						  & $10^{-4}$ & $133$ ($33287$)  & $773$ ($43459$)   & $2967$ ($45532$) & $3001$ ($17192$) & $3001$ ($5774$) \\   
    					  & $10^{-2}$ & $141$ ($36327$)  & $695$ ($41605$)   & $2982$ ($45144$) & $3001$ ($16773$) & $3001$ ($5657$) \\   
    					  & $10^0$    & $366$ ($105246$) & $1816$ ($111467$) & $3001$ ($22961$) & $3001$ ($9060$)  & $3001$ ($3623$) \\   
                          & $10^2$    & $280$ ($2117.4$)   & $1110$ ($2113$)   & $2608$ ($2014$)  & $3001$ ($2013$)  & $2979$ ($2011$) \\   
                          & $10^4$    & $65$ ($22$)      & $77$ ($22$)       & $74$ ($22$)      & $67$ ($22$)      & $64$ ($22$)     \\   
                          & $10^6$    & $12$ ($1.6$)       & $12$ ($1.6$)        & $11$ ($1.6$)       & $11$ ($1.6$)       & $10$ ($1.6$)      \\   

\bottomrule
\end{tabular}
\caption{Number of iterations (and condition number estimates) of a CG
  solver of the system~f{eq:mpt:vf} with an algebraic multigrid
  (Hypre AMG) preconditioner of the form~\eqref{eq:mpt:precond1}.}
\label{tab:mpt:2net:nontransformed}
\end{table}
\begin{table}[t]
\centering
  \begin{tabular}{l|l|l|l|l|l|l}
 $\xi_{1\leftarrow2}$     & $K_2$        & $N=8$            & $N=16$           & $N=32$           & $N=64$           & $N=128$          \\
\toprule
\multirow{7}{*}{$10^4$}   & $10^{-6}$     & $8$ ($1.2$) & $9$ ($1.2$) & $9$ ($1.2$) & $9$ ($1.2$) & $9$ ($1.2$) \\
  						  & $10^{-4}$     & $8$ ($1.2$) & $9$ ($1.2$) & $9$ ($1.2$) & $9$ ($1.2$) & $9$ ($1.2$) \\
   						  & $10^{-2}$     & $8$ ($1.2$) & $9$ ($1.2$) & $9$ ($1.2$) & $8$ ($1.2$) & $7$ ($1.1$) \\
    					  & $10^0$        & $8$ ($1.2$) & $8$ ($1.1$) & $6$ ($1.1$) & $6$ ($1.1$) & $6$ ($1.1$) \\
    					  & $10^2$        & $8$ ($1.1$) & $7$ ($1.1$) & $6$ ($1.1$) & $6$ ($1.1$) & $6$ ($1.1$) \\
    					  & $10^4$        & $7$ ($1.1$) & $6$ ($1.1$) & $6$ ($1.1$) & $6$ ($1.1$) & $7$ ($1.1$) \\
    					  & $10^6$        & $7$ ($1.1$) & $6$ ($1.1$) & $6$ ($1.1$) & $6$ ($1.1$) & $7$ ($1.1$) \\
\midrule
 \multirow{7}{*}{$10^{6}$} & $10^{-6}$     & $8$ ($1.2$) & $9$ ($1.2$) & $9$ ($1.2$) & $9$ ($1.2$) & $9$ ($1.2$) \\
  						   & $10^{-4}$     & $8$ ($1.2$) & $9$ ($1.2$) & $9$ ($1.2$) & $9$ ($1.2$) & $9$ ($1.2$) \\
  						   & $10^{-2}$     & $8$ ($1.2$) & $9$ ($1.2$) & $9$ ($1.2$) & $9$ ($1.2$) & $9$ ($1.2$) \\
  						   & $10^0$        & $8$ ($1.2$) & $9$ ($1.2$) & $9$ ($1.2$) & $9$ ($1.2$) & $8$ ($1.1$) \\
  						   & $10^2$        & $8$ ($1.2$) & $9$ ($1.2$) & $9$ ($1.2$) & $8$ ($1.2$) & $7$ ($1.1$) \\
  						   & $10^4$        & $8$ ($1.2$) & $9$ ($1.2$) & $9$ ($1.2$) & $8$ ($1.2$) & $7$ ($1.1$) \\
  						   & $10^6$        & $8$ ($1.2$) & $8$ ($1.2$) & $8$ ($1.2$) & $8$ ($1.2$) & $7$ ($1.1$) \\
\bottomrule
\end{tabular}
\caption{Number of iterations (and condition number estimates) of a CG
  solver of the system~\eqref{eq:mpt:vf:transformed} with an algebraic
  multigrid (Hypre AMG) preconditioner of the
  form~\eqref{eq:mpt:transformed:precond}.}
\label{tab:mpt:2net:transformed}
\end{table}
We first compare the performance of the preconditioners
\eqref{eq:mpt:precond1} and \eqref{eq:mpt:transformed:precond}. We let
$K_1=1.0$, and consider different values of the parameters $K_2,
\xi_{1 \leftarrow 2}$ and different mesh resolutions $N$. For the
standard formulation (Table~\ref{tab:mpt:2net:nontransformed}), the
number of iterations (and condition number) is not bounded and
increases with the ratio between  $\xi_{1 \leftarrow 2}$ and $K_2$. 
We see that the growth is somewhat less than the predicted linear growth. 
In contrast, for the new formulation (Table~\ref{tab:mpt:2net:transformed}), we observe that both the
number of iterations and the condition number stays nearly constant across the whole range of parameter values tested.
\end{example}
\begin{example}
\label{ex:mpt:3net}
In this final example, we study the performance of the preconditioner
\eqref{eq:mpt:transformed:precond} for three networks. We report the
results for $K_1=1.0$, and different values of the parameters $K_2,
K_3, \xi_{1 \leftarrow 2}, \xi_{1 \leftarrow 3}, \xi_{2 \leftarrow 3}
= (10^{-4}, 10^{-2}, 10^{0}, 10^{2},10^{4})$ and different mesh resolutions $N = (16, 32, 64)$. The results are shown in
Figure~\ref{fig:ex:3net}. We observe that the number of iterations
stays between $4$ and $6$ across the whole range of parameters tested,
with condition numbers estimated in the range $1.0-1.25$. 
\begin{figure}[h] 
  \includegraphics[width=\textwidth]{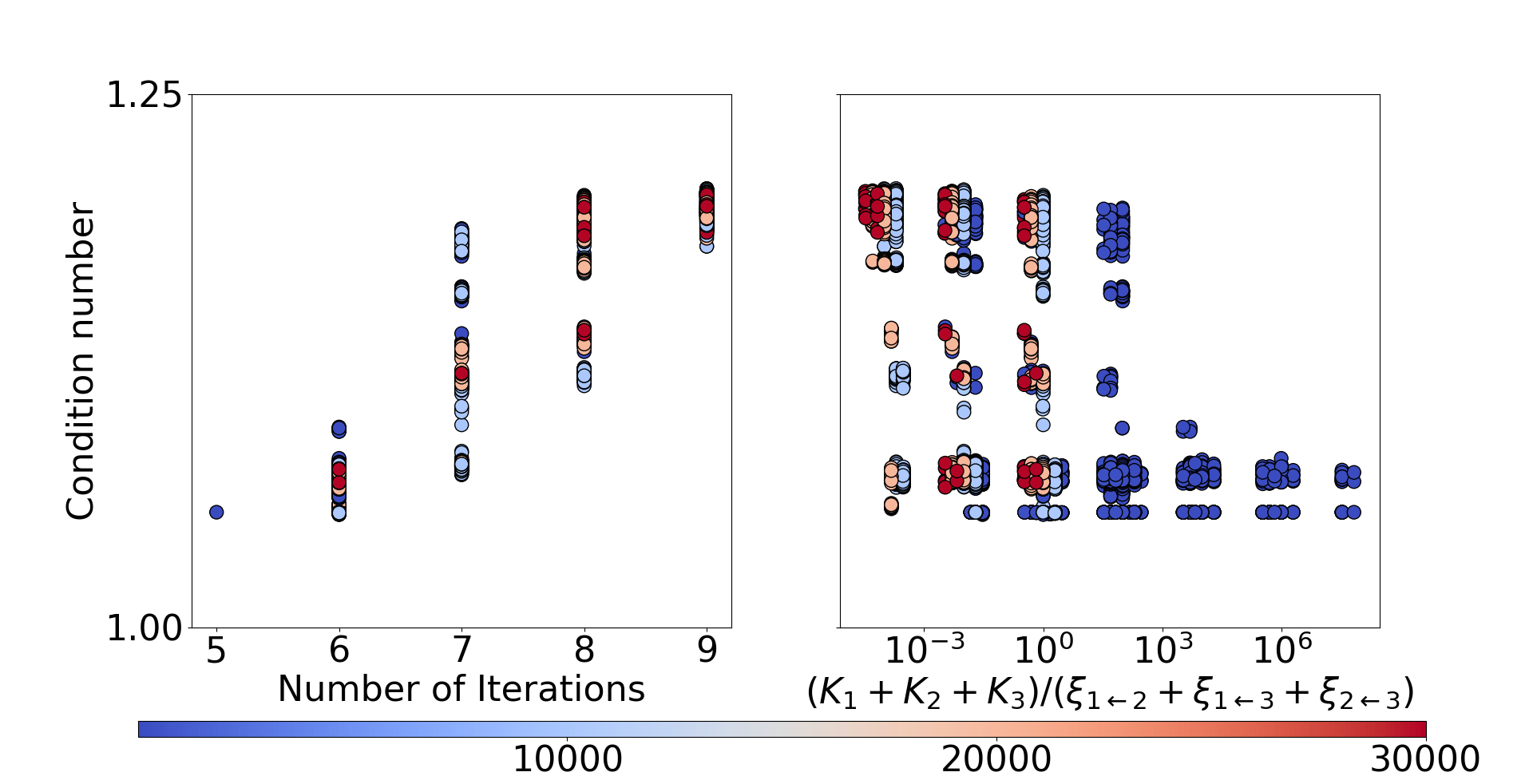}
  \caption{Example~\ref{ex:mpt:3net}: each point on the graphs
    represents a simulation performed with different parameters. The
    color represents the magnitude of $\xi_{1 \leftarrow 2}+\xi_{1
      \leftarrow 3}+\xi_{2 \leftarrow 3}$ from smaller (blue) to
    larger (red). Left: the condition number of the operator versus
    the number of iterations. Right: condition number versus the ratio
    between the sum of $\xi_{j\leftarrow i}$ and sum of $K_j$ (x-axis is logarithmic y-axis is linear). }
  \label{fig:ex:3net}
\end{figure}
\end{example}
\section{Conclusion}
In this paper we have introduced a transformation, based on the
congruence of the involved matrices, that transforms MPT systems to a
form where diagonal block preconditioners are highly effective.  The
transformation removes a problem that elliptic systems may have when
the elliptic constant is small compared to the continuity constant
because of large low order terms.

\begin{acknowledgement}
The authors would like to acknowledge Jeonghun J. Lee (Baylor
University) and Travis Thompson (University of Oxford) for
constructive discussions related to this work.
\end{acknowledgement}

\end{document}